\documentclass[12pt]{article}

\usepackage[margin=1in]{geometry}
\usepackage{amscd}
\usepackage{amssymb}
\usepackage{amsmath}
\usepackage{amsthm}
\usepackage{enumerate}
\usepackage{tikz}
\usepackage{wrapfig, framed, caption}
\usepackage{float}
\usetikzlibrary{arrows}
\usepackage[font=small,labelfont=bf]{caption}
\usepackage{xcolor}
\usepackage{mathtools}
\usepackage[colorlinks=true, linkcolor=blue, citecolor=blue,
pagebackref=true]{hyperref}
\usepackage{fouriernc}
\usepackage{ulem}

\newtheorem{theorem}{Theorem}

\newtheorem*{theoremY*}{Theorem Y}

\newtheorem*{theoremAB*}{Theorem AB}

\newtheorem*{linearformsmtp*}{Mass transference principle for linear forms}
\newtheorem{corollary}{Corollary}
\newtheorem*{corollary*}{Corollary}

\newtheorem*{claim*}{Claim}
\newtheorem{conjecture}{Conjecture}

\theoremstyle{definition}
\newtheorem{definition}{Definition}
\theoremstyle{remark}
\newtheorem{remark}{Remark}
\newtheorem*{remark*}{Remark}
\newtheorem{example}{Example}

\renewcommand{\Bbb}[1]{\mathbb{#1}}

\newcommand{\N}{{\Bbb N}}         

\newcommand{\R}{{\Bbb R}}        

\newcommand{\Z}{{\Bbb Z}}         

\newcommand{\cA}{{\cal A}}
\newcommand{\cB}{{\cal B}}

\newcommand{\cF}{{\cal F}}
\newcommand{\cG}{{\cal G}}
\newcommand{\cH}{{\cal H}}

\newcommand{\cK}{{\cal K}}
\newcommand{\cL}{{\cal L}}
\newcommand{\cM}{{\cal M}}

\newcommand{\cR}{{\cal R}}

\newcommand{\ba}{\mathbf{a}}
\newcommand{\bi}{\mathbf{i}}
\newcommand{\bj}{\mathbf{j}}
\newcommand{\x}{\mathbf{x}}
\newcommand{\y}{\mathbf{y}}
\newcommand{\p}{\mathbf{p}}

\newcommand{\bt}{\mathbf{t}}

\newcommand{\diam}{\textrm{diam}}
\newcommand{\dist}{\operatorname{dist}}

\DeclareMathOperator{\dimh}{\dim_H}
\DeclareMathOperator{\dimS}{\dim_S}
\DeclareMathOperator{\supp}{supp}

\title{A survey of recent extensions and generalisations of the Mass Transference Principle}

\author{Demi Allen \\ (University of Exeter) \and Edouard Daviaud \\ (Université Paris-Est Créteil)}

\date{}

\begin{document}
\frenchspacing
\maketitle

\begin{abstract}
In this article we survey some recent extensions and generalisations of the celebrated Mass Transference Principle (Beresnevich and Velani, Annals of Mathematics, 2006).  
\end{abstract}

{\footnotesize
\tableofcontents
}

\section{Introduction}

The Mass Transference Principle established by Beresnevich and Velani in 2006 \cite{BV MTP} is a powerful result enabling one to pass between Lebesgue measure and Hausdorff measure statements for related $\limsup$ sets of balls. This is surprising since the latter, Hausdorff measure, is generally viewed as being a refinement of the former, Lebesgue measure. Indeed, for subsets of $\R^k$, the $k$-dimensional Hausdorff measure is a constant multiple times the \mbox{$k$-dimensional} Lebesgue measure \cite{Falconer}. The Mass Transference Principle has cemented itself as a fundamental result in Metric Number Theory which has been found to have wide-ranging applications. Moreover, the Mass Transference Principle has now been extended in a number of directions and further generalisations and applications are constantly emerging. Indeed, only a few years ago, a survey was written on generalisations and applications of the Mass Transference Principle up until that point \cite{Allen-Troscheit} but, already, this is somewhat outdated. In this survey, we will focus on collecting together various generalisations of the Mass Transference Principle, paying particular attention to the advances made in the last few years since the writing of \cite{Allen-Troscheit}. There have also been many newfound applications of the various mass transference principles we discuss here. However, due to space constraints, we will not delve into any of these in any detail in this particular survey.

\section{Preliminaries}

In this section, we will establish some notation and gather together some of the necessary preliminaries that will be used throughout. The typical formulation of a mass transference principle result is that it allows us to infer information regarding (usually) Hausdorff measures of limsup sets in some setting from \emph{a priori} weaker statements regarding (usually) Lebesgue measure of some related limsup set. We recall that if we have a sequence, say $(A_i)_{i \in \N}$, of subsets of some metric space $X$, then the \emph{limsup} of this sequence is
\begin{align*}
\limsup_{i \to \infty}{A_i} &= \{x \in X: x \in A_i \text{ for i.m. } i \in \N\} \\
                            &= \bigcap_{N=1}^{\infty}{\bigcup_{i=N}^{\infty}{A_i}}.
\end{align*}
Throughout, the shorthand ``i.m.'' stands for ``infinitely many''.

\subsection{Hausdorff measures and dimension} 

We define a \emph{dimension function} to be any continuous function $f: \R_{> 0} \to \R_{> 0}$ such that \mbox{$f(r) \to 0$ as $r \to 0$}. Recall that for a function $f: \R_{>0} \to \R_{>0}$, we say that $f$ is \emph{doubling} if there exists a constant $\lambda > 1$ such that for all $x > 0$ we have 
\[f(2x) \leq \lambda f(x).\]

For a subset $F$ of $\R^k$ and a real number $\rho > 0$, we say that a countable collection of balls $\{B_i = B(x_i, r_i)\}_{i \in I}$ is a \emph{$\rho$-cover} for $F$ if $r_i < \rho$ for every $i \in I$ and $F \subset \bigcup_{i \in I}{B_i}$. For a dimension function $f$ and a real number $\rho > 0$, the \emph{$\rho$-approximate Hausdorff $f$-measure} of $F$ is defined as 
\[\cH_{\rho}^{f}(F) = \inf\left\{\sum_{i}{f(r_i)}: \{B_i\} \text{ is a $\rho$-cover of } F\right\}.\]
The \emph{Hausdorff $f$-measure} of $F$ is then defined as 
\[\cH^f(F) = \lim_{\rho \to 0}{\cH_{\rho}^{f}{(F)}}.\]
Note that this limit exists (although may be infinite) since as $\rho$ decreases, the number of possible $\rho$-covers is restricted and so the value of $\cH_{\rho}^{f}(F)$ is monotonically increasing as $\rho$ decreases. When $f(r) = r^s$ for some $s \geq 0$, we write $\cH^s$ in place of $\cH^{f}$. The Hausdorff dimension of $F$ is
\[\dimh{F} = \inf\{s \geq 0: \cH^s(F) = 0\}.\]

\subsection{Dimension of measures}

Let $\cB(\R^k)$ denote the set of all Borel subsets of $\R^k$ and let $\cM(\R^k)$ denote the set of all Borel probability measures on $\R^k$. For a measure $\mu \in \cM(\R^k)$, we write $\supp(\mu)$ to denote the \emph{support} of $\mu$.
 
Suppose $\mu \in \cM(\R^k)$ and let $x\in \mbox{supp}(\mu)$. Then the \emph{lower} and \emph{upper local dimensions} of $\mu$ at $x$ are  defined, respectively, as
\begin{align*}
\underline\dim(\mu,x)&=\liminf_{r\to 
 0}\frac{\log(\mu(B(x,r)))}{\log(r)} \; , \\[2ex]
 \overline\dim(\mu,x)&=\limsup_{r\to 
 0}\frac{\log (\mu(B(x,r)))}{\log(r)}.
 \end{align*}
The \emph{lower Hausdorff dimension} of the measure $\mu$ is given by
\begin{equation*}
\underline{\dim}(\mu)=\inf\left\{\dimh(E):\, E\in\mathcal{B}(\mathbb{R}^d),\, \mu(E)>0\right\}. 
\end{equation*}
If $\underline{\dim}(\mu,x)=\overline{\dim}(\mu,x)$ for $\mu$-almost every $x \in\mbox{supp}(\mu)$, then this common value is simply denoted by $\dim(\mu)$ and $\mu$ is said to be \emph{exact dimensional}. For further details see, for example, \cite{Falconer} or \cite{Falconer_techniques}. 

More generally, at various points, we will also consider probability measures on locally compact metric spaces. If $X$ is a locally compact metric space equipped with a probability measure $\mu$, we say that $\mu$ is \emph{$\delta$-Ahlfors-David regular} if there exist constants $c_1,c_2,r_0>0$ such that for any $x \in X$ and for any $0 < r < r_0$, we have 
\[c_1 r^{\delta} \leq \mu(B(x,r)) \leq c_2 r^{\delta}.\]
If a measure $\mu$ is $\delta$-Ahlfors-David regular for \emph{some} $\delta$, we may simply say that $\mu$ is \emph{Ahlfors-David regular} if the particular value of $\delta$ is not of consequence.

\subsection{Multifractal analysis of measures}
The multifractal analysis of measures consists of associating a spectrum that is correlated to the distribution of the possible local dimensions of the measure in order to finely characterise its behaviour. More precisely, given $\mu\in\cM(\R^k)$ and a real number $h>0$, it usually relies on estimating the Hausdorff dimension of the following sets:
\begin{align*}
E_{\mu}(h) &=\left\{x\in \mbox{supp}(\mu):\underline{\dim}(\mu,x)=h\right\},\nonumber \\[2ex]
\widetilde{E}_{\mu}(h) &=\left\{x\in \mbox{supp}(\mu):\dim(\mu,x)=h\right\}.
\end{align*}

Let $\mu \in\mathcal{M}(\mathbb{R}^d)$ and define
\[\Theta_{\mu}(q,r)= \inf\left\{\sum_{i\in I} \mu(B(x_i,r))^q\right\}.\]
where the infimum is taken over all countable collections of balls $\left\{B(x_i,r) \right\}_{i \in I}$ satisfying the following two properties: 
\begin{enumerate}[(i)]
    \item $x_i \in \supp(\mu)$ for every $i \in I$, and
    \item $B(x_i,r) \cap B(x_j,r) =\emptyset$ for $i\neq j.$
\end{enumerate}

We then define 
 \[\tau_{\mu}(q)=\liminf_{r\to 0}\frac{\log (\Theta_{\mu}(q,r))}{ \log r}.\]
One often aims to correlate the dimension of the above sets with the scaling functions $\tau_{\mu}(q)$.

We say that the measure $\mu $ satisfies the \emph{multifractal formalism} at $h$ when  $$\dimh(E_{\mu}(h))=\tau_{\mu}^* (h),$$ where $\tau_{\mu} ^{*}(h):=\inf_{q \in \R}\left\{hq-\tau_{\mu} (q)\right\}$ denotes the Legendre transform of $\tau_{\mu}$ at $h$.

\subsection{Iterated function systems and self-similar measures}
Let $D$ be a closed subset of $\R^k$. A function $\varphi:D \to D$ is said to be a \emph{contraction} if there exists some constant $0 < r < 1$ such that 
\[|\varphi(x)-\varphi(y)| \leq r|x-y|\]
for all $x,y \in D$. If we have equality above then we say that $\varphi$ is a \emph{contracting similarity} and that $r$ is the \emph{contraction ratio} of $\varphi$. We call a finite collection $\Phi = \{\varphi_i:D \to D\}_{i=1}^{n}$ of at least two contractions an \emph{iterated function system} or \emph{IFS} for short. If all of the maps in the iterated function system $\Phi$ are contracting similarities, we say that $\Phi$ is a \emph{self-similar IFS}. It is a well-known result due to Hutchinson \cite{Hutchinson} that, for any iterated function system $\Phi = \{\varphi_i:D \to D\}_{i=1}^{n}$, there exists a unique, non-empty, compact set $X \subset D$ such that
\[X = \bigcup_{i=1}^{n}{\varphi_i(X)}.\]
We call $X$ the \emph{attractor} of the iterated function system $\Phi$. If $\Phi$ is a self-similar IFS, then $X$ is a \emph{self-similar set}. In this case, we may still refer to the attractor of $\Phi$, but we may alternatively refer to $X$ as being the \emph{self-similar set associated to $\Phi$}. Moreover, given any real numbers $0<p_1,...,p_n <1$ satisfying $p_1+...+p_n=1,$ it is also known (\cite{Hutchinson}) that there exists a unique measure $\mu\in\mathcal{M}(\mathbb{R}^k)$ with $\supp{\mu} = X$ satisfying
$$\mu(\cdot)=\sum_{i=1}^n p_i \mu(\varphi_i^{-1}(\cdot)).$$
We call such measures \emph{self-similar measures} and, due to a result of Feng and Hu \cite{FH}, we know that these measures are always exact-dimensional.

\subsection{The open set condition}

Suppose that $\Phi = \{\varphi_i: D \to D\}_{i=1}^{n}$ is an iterated function system. Let us write $\Lambda = \{1,\dots,n\}$ and define
\[\Lambda^*= \bigcup_{N=0}^{\infty}{\Lambda^N}.\]
Thus, $\Lambda^*$ consists of all possible finite sequences consisting of entries from $\Lambda$. For $\bi = (i_1,\dots,i_{|\bi|}) \in \Lambda^*$, we adopt the shorthand notation
\[\varphi_{\bi} = \varphi_{i_1} \circ \varphi_{i_2} \circ \dots \circ \varphi_{i_{|\bi|}}.\] 

We say that $\Phi$ satisfies the \emph{open set condition} if there exists a non-empty bounded open set $U \subset D$ such that
\[\varphi_i(U) \subset U \quad \text{for each } 1 \leq i \leq n \quad \text{and} \quad \varphi_i(U) \cap \varphi_j(U) = \emptyset \quad \text{for } i \neq j.\]

\subsection{Hausdorff dimension of self-similar sets}

Suppose that $\Phi = \{\varphi_i\}_{i=1}^{n}$ is a self-similar IFS and that, for each $1 \leq i \leq n$, the map $\varphi_i$ has contraction ratio $r_i$. The \emph{similarity dimension} of $\Phi$, denoted by $\dimS{\Phi}$, is defined to be
the unique real number satisfying 
\begin{equation}
   \label{simidim}
    \sum_{i=1}^n r_i^{\dimS{\Phi}}=1.
\end{equation}

In many natural situations, the Hausdorff dimension of a self-similar set turns out to be equal to the similarity dimension of the corresponding iterated function system. In particular, it is well-known (see, for example, \cite{Falconer}) that if $\Phi = \{\varphi_i\}_{i=1}^{n}$ is a self-similar IFS satisfying the open set condition and $X$ is the associated self-similar set, then
\begin{align*} \label{similarity dimension equivalences}
\dimh{X} = \dimS{\Phi}.
\end{align*}

\subsection{Singular Value Function} \label{singular value function section}
A quantity which arises naturally when estimating the Hausdorff dimension of fractal sets that are defined by rectangles (such as self-affine sets or limsup sets generated by rectangles) is the so-called \emph{singular value function}. Here we follow the definition as given in \cite{Falconer} but specialise to the case where the sets of interest are rectangles. Suppose $\x = (x_1,\dots,x_k) \in \R^k$, $r>0$, and $1=\tau_1 \leq ... \leq \tau_k$. Define
\[R=\prod_{i=1}^k{B(x_i,r^{\tau_i})}\] 
to be the rectangle centered at $\x$ with side lengths $r^{\tau_1},...,r^{\tau_k}$. For $s \geq 0$, the $s$-dimensional \emph{singular value function} associated with $R$ is defined as
\begin{align}
\label{definitionsingularvalue}
\phi_s (R)&=r^{\max_{1\leq i\leq k}s\tau_i -\sum_{j=1}^i \tau_i -\tau_j}
\\[1ex]
&=r^{s\tau_m -m\tau_m +\sum_{j=1}^{m}\tau_j},\nonumber
\end{align}
where $m$ is the integer satisfying $m<s\leq m+1.$

 Given a dimension function $f: \R_{>0} \to \R_{>0}$, the $f$-dimensional Hausdorff content is denoted by $\cH^f_{\infty}(E)$ and defined as
$$\cH^f_{\infty}(E)=\inf\left\{\sum_{n \in \N}{f(r_n)}: (B(x_n,r_n))_{n\in\mathbb{N}} \text{ is a sequence of closed 
 balls such that } E\subset \bigcup_{n \in \N}B(x_n,r_n)\right\}.$$
When $f(r)=r^s$, $\cH^f_{\infty}$ is denoted $\cH^s_{\infty}.$ 

As is pointed out by Koivusalo and Rams in  \cite{KR}, up to some multiplicative constant, the $s$-dimensional singular value of a rectangle $R$ is nothing more than the $s$-dimensional \emph{Hausdorff content} of $R$.

\section{The ``original'' Mass Transference Principle and the Duffin--Schaeffer Conjecture} \label{MTP section}

\subsection{The Duffin--Schaeffer Conjecture}
The Mass Transference Principle proved by Beresnevich and Velani in \cite{BV MTP} was borne out of a desire to find a Hausdorff measure analogue of the famous Duffin--Schaeffer Conjecture in Diophantine Approximation. Given a function $\psi: \N \to \R_{\geq 0}$, which we will often refer to as an \emph{approximating function}, the \emph{$\psi$-well approximable points} in $\R$ are defined to be
\begin{align*}
\cA(\psi) = \left\{x \in [0,1]: \left|x-\frac{p}{q}\right| < \frac{\psi(q)}{q} \text{ for i.m. } (p,q) \in \Z \times \N \right\}.    
\end{align*}

A classical result due to Khintchine \cite{Khintchine1924} (subsequently modernised by others --- see, for example, \cite{BRV}) gives us the following characterisation of the Lebesgue measure of sets of $\psi$-well approximable points. We write $\cL(X)$ to denote the Lebesgue measure of a subset $X \subset \R^k$.

\begin{theorem}[Khintchine's Theorem]\label{Khintchine} 
For $\psi: \N \to \R_{\geq 0}$, we have
\begin{equation*}
\cL(\cA(\psi)) =
    \begin{cases}
      0 &\textrm{if } \sum_{q=1}^{\infty}{\psi(q)} < \infty \;, \\[2ex]
      1 &\textrm{if } \sum_{q=1}^{\infty}{\psi(q)} =\infty \quad \textrm{and $\psi$ is monotonic.}
    \end{cases}
\end{equation*}
\end{theorem}

The convergence part of Khintchine's Theorem can be proved easily via a straightforward application of the First Borel--Cantelli Lemma (for a statement of the First Borel--Cantelli Lemma, we refer the reader to \cite[Lemma 1.2]{Harman}). In particular, note that the proof of the convergence part of Khintchine's Theorem requires no monotonicity assumptions on the approximating function $\psi$. On the other hand, it was shown by Duffin and Schaeffer \cite{Duffin-Schaeffer} that the monotonicity assumption cannot be removed from the divergence part of Khintchine's Theorem. In the same paper where they exhibited a counterexample demonstrating this, they also formulated a conjecture for what ought to be true instead of Khintchine's Theorem when we do not assume that the approximating function $\psi$ is monotonic. To this end, given $\psi: \N \to \R_{\geq 0}$, we define 
\begin{align*}
\cA'(\psi) = \left\{x \in [0,1]: \left|x-\frac{p}{q}\right| < \frac{\psi(q)}{q} \text{ for i.m. } (p,q) \in \Z \times \N \text{ with } \gcd(p,q) = 1 \right\}.
\end{align*}

Let us denote by $\varphi$ the \emph{Euler totient function}; that is, for $n \in \N$, 
\[\varphi(n) = \#\{1 \leq k \leq n: \gcd(k,n) = 1\}.\] 
Duffin and Schaeffer predicted the following \cite{Duffin-Schaeffer}:
\begin{conjecture}[Duffin--Schaeffer Conjecture] \label{DS Conjecture} For any approximating function $\psi: \N \to \R_{\geq 0}$, we have 
\begin{equation*}
\cL(\cA'(\psi)) =
    \begin{cases}
      0 &\textrm{if } \sum_{q=1}^{\infty}{\frac{\varphi(q)\psi(q)}{q}} < \infty \;, \\[2ex]
      1 &\textrm{if } \sum_{q=1}^{\infty}{\frac{\varphi(q)\psi(q)}{q}} =\infty.
    \end{cases}
\end{equation*}
\end{conjecture}

In fact, the convergence part of the above conjecture follows easily from the First Borel--Cantelli Lemma. Thus, it is the divergence statement which is the real substance of the above conjecture. This eluded the mathematical community for almost 80 years, until there was a recent major breakthrough and the conjecture was proved by Koukoulopoulos and Maynard \cite{KM}.

\begin{theorem} \label{DS theorem}
The Duffin--Schaeffer Conjecture is true.    
\end{theorem}

In the intervening years between the formulation of Conjecture \ref{DS Conjecture} by Duffin and Schaeffer \cite{Duffin-Schaeffer} and its eventual resolution by Koukoulopoulos and Maynard \cite{KM}, a natural higher-dimensional analogue of the Duffin--Schaeffer Conjecture was posed by Sprind\v{z}uk and proved by Pollington and Vaughan \cite{PV1990}.

For $k \in \N$ and $\psi: \N \to \R_{\geq 0}$, we define the \emph{simultaneously $\psi$-well approximable points} in $\R^k$ as:
\begin{align*}
\cA_{k}(\psi) = \left\{x \in [0,1]^k: \|qx-\p\| < \psi(q) \text{ for i.m. } (\p,q) \in \Z^k \times \N \right\},    
\end{align*}
where $\|\cdot\|$ denotes the supremum norm. In this case, a characterisation of the Lebesgue measure of sets of simultaneously $\psi$-well approximable points is again given by a result due to Khintchine \cite{Khintchine1925} (and subsequently modernised, see \cite{BRV} for details).

\begin{theorem}[Khintchine's Theorem in $\R^k$] \label{higher dimensional Khintchine}
For $\psi: \N \to \R_{\geq 0}$, we have
\begin{equation*}
\cL(\cA_{k}(\psi)) =
    \begin{cases}
      0 &\textrm{if } \sum_{q=1}^{\infty}{\psi(q)^k} < \infty \;, \\[2ex]
      1 &\textrm{if } \sum_{q=1}^{\infty}{\psi(q)^k} =\infty \quad \textrm{and $\psi$ is monotonic.}
    \end{cases}
\end{equation*}
\end{theorem}

\begin{remark}
It follows from a result of Gallagher \cite{Gallagher} that monotonicity of the approximating function $\psi$ is not actually required in Theorem \ref{higher dimensional Khintchine} when $k \geq 2$.
\end{remark}

To state the higher-dimensional version of the Duffin--Schaeffer Conjecture, for $k \in \N$ and given $\psi: \N \to \R_{\geq 0}$, we define
\begin{align*}
\cA_{k}''(\psi) &= \left\{x \in [0,1]: \|qx-\p\| < \psi(q) \text{ for i.m. } (\p,q)=(p_1,\dots,p_k,q) \in \Z^k \times \N\right. \\
&\phantom{====================}\left.\text{ with } \gcd(p_i,q) = 1 \text{ for each } 1 \leq i \leq k \right\}.
\end{align*}

Sprind\v{z}uk postulated the following conjecture regarding the Lebesgue measure of the sets $\cA_{m}''(\psi)$ when $\psi$ is not assumed to be monotonic. 

\begin{conjecture} \label{higher-dimensional DSC} For any approximating function $\psi: \N \to \R_{\geq 0}$, we have
\begin{equation*}
\cL(\cA_{k}''(\psi)) =
    \begin{cases}
      0 &\textrm{if } \sum_{q=1}^{\infty}{\left(\frac{\varphi(q)\psi(q)}{q}\right)^k} < \infty \;, \\[2ex]
      1 &\textrm{if } \sum_{q=1}^{\infty}{\left(\frac{\varphi(q)\psi(q)}{q}\right)^k} =\infty.
    \end{cases}
\end{equation*}
\end{conjecture}

Notice that when $k=1$, Conjecture \ref{higher-dimensional DSC} matches Conjecture \ref{DS Conjecture} and so in this sense is a natural generalisation of it. In the cases when $k \geq 2$, Conjecture \ref{higher-dimensional DSC} has been verified by Pollington and Vaughan \cite{PV1990}. Of course, putting this together now with Theorem \ref{DS theorem} gives a full resolution of Conjecture \ref{higher-dimensional DSC}.

\begin{theorem} \label{higher dimensional DS theorem}
Conjecture \ref{higher-dimensional DSC} is true. 
\end{theorem}

\subsection{The ``original'' Mass Transference Principle and the Hausdorff measure Duffin--Schaeffer Conjecture}

In \cite{BV MTP}, Beresnevich and Velani were aiming to present natural Hausdorff measure generalisations of Conjectures \ref{DS Conjecture} and \ref{higher-dimensional DSC}. The natural, all-encompassing, statement presented in \cite{BV MTP} was the following:

\begin{conjecture}\label{Hausdorff DSC}
Let $k \in \N$ and let $f$ be a dimension function such that $r^{-k}f(r)$ is monotonic. For any approximating function $\psi: \N \to \R_{\geq 0}$, we have
\begin{equation*}
\cH^f(\cA_{k}''(\psi)) =
    \begin{cases}
      0 &\textrm{if } \sum_{q=1}^{\infty}{f\left(\frac{\psi(q)}{q}\right)\varphi(q)^k} < \infty \;, \\[2ex]
      \cH^f([0,1]^k) &\textrm{if } \sum_{q=1}^{\infty}{f\left(\frac{\psi(q)}{q}\right)\varphi(q)^k} =\infty.
    \end{cases}
\end{equation*}
\end{conjecture}

We remark that the convergence part of Conjecture \ref{Hausdorff DSC} can be proved via a relatively straightforward covering argument and so the real essence of the conjecture is the divergence part. Moreover, when $\psi$ is monotonic, Conjecture \ref{Hausdorff DSC} is a classical result due to Jarn\'{\i}k~\cite{Jarnik1931}. Thus, the really interesting case is when $\psi$ is not assumed to be monotonic.

In the case when $f(r)=r^k$, i.e. essentially Lebesgue measure, Conjecture \ref{Hausdorff DSC} is equivalent to Conjecture \ref{higher-dimensional DSC}. Therefore, on the face of it, it would seem that Conjecture \ref{Hausdorff DSC} really is more general than Conjectures \ref{DS Conjecture} and \ref{higher-dimensional DSC}. Moreover, a proof of Conjecture \ref{Hausdorff DSC} would immediately imply the truth of the other two conjectures. Rather surprisingly, the Mass Transference Principle shows that, in fact, the converse of this is also true; that a proof of Conjecture \ref{higher-dimensional DSC} would imply the seemingly more general Conjecture \ref{Hausdorff DSC}.

Given $k \in \N$, a dimension function $f: \R_{> 0} \to \R_{> 0}$, and a ball $B = B(x,r)$ in $\R^k$, we define an associated ball
\[B^f = B(x,f(r))^{\frac{1}{k}}.\]
When $f(r) = r^s$ for some $s \geq 0$, we write $B^s = B^f$. Also note that with this notational convention, we have $B^k = B$. The following Mass Transference Principle was established by Beresnevich and Velani in \cite[Theorem 2]{BV MTP}.

\begin{theorem}[The Mass Transference Principle] \label{MTP}
Let $\{B_i = B(x_i, r_i)\}_{i \in \N}$ be a sequence of balls in $\R^k$ with $r_i \to 0$ as $i \to \infty$. Let $f: \R_{> 0} \to \R_{> 0}$ be a dimension function such that $r^{-k}f(r)$ is monotonic and suppose that for any ball $B$ in $\R^k$, we have 
\[\cH^k\left(B \cap \limsup_{i \to \infty}{B_i^f}\right) = \cH^k(B).\]
Then, for any ball $B$ in $\R^k$, we have
\[\cH^f\left(B \cap \limsup_{i \to \infty}{B_i}\right) = \cH^f(B).\]
\end{theorem}

It was shown already in \cite{BV MTP} that a consequence of Theorem \ref{MTP} was that a proof of Conjecture \ref{higher-dimensional DSC} would imply Conjecture \ref{Hausdorff DSC}. Since we now have a full resolution of Conjecture \ref{higher-dimensional DSC}, with the cases when $k \geq 2$ due to Pollington and Vaughan \cite{PV1990}, and the case when $k=1$ due to Koukoulopoulos and Maynard \cite{KM}, we also have:

\begin{theorem} \label{Hausdorff measure DS theorem}
Conjecture \ref{Hausdorff DSC} is true.
\end{theorem}

We refer the reader to, for example, \cite{BV MTP} for further details on how Theorem \ref{Hausdorff measure DS theorem} follows from Theorem \ref{higher dimensional DS theorem}. Similar arguments also yield further surprising consequences of the Mass Transference Principle. For example, one can deduce Jarn\'{\i}k's Theorem \cite{Jarnik1931} from Khintchine's Theorem and the Jarn\'{\i}k--Besicovitch Theorem (\cite{Jarnik1929}, \cite{Besicovitch}) can be inferred from Dirichlet's Theorem. We will not spend time elaborating on these particular applications here as our primary focus is to survey some of the more recent developments surrounding the Mass Transference Principle, but we instead refer the interested reader to the likes of \cite{BV MTP} and \cite{BRV} for further details.

\subsection{A general mass transference principle from ``balls to balls''}

We conclude our section on the ``original'' Mass Transference Principle by briefly mentioning a more general version of Theorem \ref{MTP} also recorded by Beresnevich and Velani in \cite{BV MTP}. Let $(X,d)$ be a locally compact metric space and let $g: \R_{> 0} \to \R_{> 0}$ be a doubling dimension function. Further suppose that there exist constants $0 < c_1 < c_2 < \infty$ and $r_0 > 0$ such that for any ball $B=B(x,r)$ with $x \in X$ and $r >0$, we have
\begin{align} \label{g doubling condition}
c_1 g(r) \leq \cH^g(B) \leq c_2 g(r).
\end{align}

Given a dimension function $f: \R_{> 0} \to \R_{> 0}$ and a ball $B(x,r)$ in $X$, define
\[B^f = B(x,g^{-1}f(r)).\]
Note that $B^g = B$.

In this more general setting, Beresnevich and Velani state the following general version of the Mass Transference Principle in \cite[Theorem 3]{BV MTP}.

\begin{theorem} \label{BV general MTP}
Let $(X,d)$ and $g$ be as above. Let $\{B_i = B(x_i, r_i)\}_{i \in \N}$ be a sequence of balls in $X$ with $r_i \to 0$ as $i \to \infty$. Let $f: \R_{> 0} \to \R_{> 0}$ be a dimension function such that $\frac{f(r)}{g(r)}$ is monotonic. Suppose that for any ball $B$ in $X$, we have
\[\cH^g\left(B \cap \limsup_{i \to \infty}{B_i^f}\right) = \cH^g(B).\]
Then, for any ball $B$ in $X$, we have 
\[\cH^f\left(B \cap \limsup_{i \to \infty}{B_i^g}\right) = \cH^f(B)\]
\end{theorem}

In the case that $X=\R^k$ and $g(x) = x^k$, Theorem \ref{BV general MTP} reduces to Theorem \ref{MTP}. The more general version given in Theorem \ref{BV general MTP}, though, has found applications in, for example, Diophantine Approximation on fractals, where Theorem \ref{MTP} cannot be applied directly. See, for example, \cite{LSV} and \cite{WWX2017}.

\section{Precursors to the Mass Transference Principle}

The idea of studying points that are approximable at certain rates using geometric properties of a ``good'' measure takes its roots in the 1990's with work of Dodson, Rynne and Vickers \cite{DRV} and work of Dodson, Melian, Pestana and Velani \cite{DMPV}. More precisely, the notion of ubiquitous systems with respect to a radius function and a measure was introduced in \cite{DRV} and refined in \cite{DMPV}. Subsequently, the theory of ubiquitous systems has been further developed and this theory has become an important topic in the context of Diophantine approximation. We refer the reader to \cite{BDV} and references therein for further information. For the moment, we will focus on the definition of a ubiquitous system as presented in \cite{DMPV}.

\begin{definition}
Fix a subset $K \subset \mathbb{R}^k$ and let $\mu$ be an Ahlfors-David regular measure on $K$. Let $\cR=\{x_{\alpha}\}_{\alpha\in J}$ be a countable collection of points in $K$ and let $\beta:J\to \R_{>0}$. For each $N\in\mathbb{N}$ define 
\[J(N)=\left\{\alpha : \ \beta(\alpha)\leq N\right\}.\] 
Suppose that there exists a positive decreasing function $\rho:\R_{> 0} \to \R_{> 0}$ such that
\begin{equation}
\label{equaubi}
\mu\left(K \setminus \bigcup_{\alpha\in J(N)}\frac{1}{3}B(x_{\alpha},\rho(N))\right)\to0 \quad \text{as } N \to \infty.
\end{equation}
Then, the system $(\mathcal{R},\beta)$ is said to be \emph{$\mu$-ubiquitous with respect to $\rho$}.
\end{definition}

In \cite{DMPV}, the following ubiquity theorem was obtained.

\begin{theorem}
Let $K$ be a subset of $\mathbb{R}^k$ and let $\mu$ be a $\delta$-Ahlfors-David regular measure on $K$. Suppose $(\mathcal{R},\beta)$ is a $\mu$-ubiquitous system with respect to a function $\rho: \R_{> 0} \to \R_{> 0}$. Then, for every $t\geq 1$,
$$\dimh \left\{y:y\in B(x_{\alpha},\rho(\beta(\alpha))^{t})\text{ for infinitely many  }\alpha \in J\right\}\geq \frac{\delta}{t}.$$
\end{theorem}
Note that in \eqref{equaubi}, the balls all have the same radius $\rho(N).$ In particular, when $\mathcal{R}$ is a sequence of points $(x_n)_{n\in \mathbb{N}}$ and $(\rho(n))_{n\in\mathbb{N}}$ 
 is a sequence of positive numbers satisfying $\rho(n) \to 0$, defining $\beta:\mathbb{N}\to\mathbb{N}$ by $\beta(n)=n$ and asking $(\mathcal{R},\beta)$ to be $\mu$-ubiquitous with respect to $(\rho(n))_{n\in\mathbb{N}}$ is \emph{a priori} a stronger requirement than asking that 
 \begin{equation}
 \label{fullmeas}
\mu\left(\limsup_{n \to \infty}B(x_n,\rho(n))\right)=1.  
 \end{equation}
However, it was only in 2000, with the work of Jaffard \cite{Jaffard1,Ja}, that the above results were refined to obtain results assuming only \eqref{fullmeas}. More explicitly, Jaffard \cite{Jaffard1,Ja} showed the following.

\begin{theorem}
Let $K$ be a subset of $\R^k$, let $\delta \geq 0$, and suppose that $\mu$ is a $\delta$-Ahlfors-David regular measure on $K$. Let $(B_n:=B(x_n,r_n))_{n\in\N}$ be a sequence of balls in $K$ such that $r_n \to 0$ as $n \to \infty$. If
\[\mu\left(\limsup_{n\to\infty}B_n\right)=1,\]
then for every $t\geq 1$, there exists a dimension function $f:\R_{> 0} \to \R_{>0}$ such that 
\[\lim_{r\to 0}\frac{\log f(r)}{\log r }=\frac{\delta}{t} \quad \text{and} \quad \cH^{f}\left(\limsup_{n\to\infty} B(x_n ,r_{n}
^{t})\right)>0.\]
In particular 
\[\dimh \left(\limsup_{n\to\infty} B(x_n ,r_n ^{t})\right) \geq \frac{\delta}{t}.\] 
\end{theorem}

\section{Mass Transference Principle for Systems of Linear Forms}

Following their discovery of the original Mass Transference Principle (Theorem \ref{MTP}), Beresnevich and Velani subsequently extended the idea and proved a mass transference principle for systems of linear forms \cite{BV slicing} by combining Theorem \ref{MTP} with a ``slicing'' result from geometric measure theory. The proof method employed resulted in some undesirable technical restrictions being imposed on the result obtained by Beresnevich and Velani in this direction (see Remark \ref{slicing technical conditions}). These technical restrictions were conjectured to be unnecessary in \cite[Conjecture E]{BBDV} and were later removed in \cite{Allen-Beresnevich}. We state here the mass transference principle for systems of linear forms given in \cite{Allen-Beresnevich}.

Let $k,m \geq 1$ and $\ell \geq 0$ be integers such that $k=m+\ell$. Let $\cR = (R_j)_{j \in \N}$ be a sequence of $\ell$-dimensional planes in $\R^k$. Let $\|\cdot\|$ be any fixed norm on $\R^k$ and, for $j \in \N$ and $\delta \geq 0$, define
\[\Delta(R_j,\delta) = \{\x \in \R^k: \dist(\x,R_j)<\delta\},\]
where 
\[\dist(\x,R_j) = \inf\{\|\x-\y\|:\y \in R_j\}.\]
Let $\Upsilon = (\Upsilon_j)_{j \in \N}$ be a sequence of non-negative real numbers such that $\Upsilon_j \to 0$ as $j \to \infty$. Consider the set
\[\Lambda(\Upsilon) = \left\{\x \in \R^k: \x \in \Delta(R_j,\Upsilon_j) \text{ for infinitely many } j \in \N\right\}.\]

\begin{theorem} \label{mtp for linear forms}
Let $\cR$ and $\Upsilon$ be as defined above. Suppose that $f: \R_{> 0} \to \R_{> 0}$ is a dimension function such that $r^{-k}f(r)$ is monotonic and the function $g: \R_{> 0} \to \R_{> 0}$ defined by $g(r) = r^{-\ell}f(r)$ is also a dimension function. Let $\Omega$ be a ball in $\R^k$ and suppose that for any ball $B$ in $\Omega$, we have
\[\cH^k\left(B \cap \Lambda\left(g(\Upsilon)^{\frac{1}{m}}\right)\right) = \cH^k(B).\]
Then, for any ball $B$ in $\Omega$, we have
\[\cH^f(B \cap \Lambda(\Upsilon)) = \cH^f(B).\]
\end{theorem}

\begin{remark} \label{slicing technical conditions}
In \cite[Theorem 4.1]{BV slicing}, Beresnevich and Velani draw the same conclusion as in Theorem \ref{mtp for linear forms} but subject to the additional condition that there exists an $m$-dimensional linear subspace $V$ of $\R^k$ satisfying:
\begin{enumerate}[(i)]
\item{$V \cap R_j \neq \emptyset$ for all $j \in \N$, and}
\item{$\sup_{j \in \N}\diam(V \cap \Delta(R_j,1))<\infty$.}
\end{enumerate}
\end{remark}

Initially, the extra conditions required by Beresnevich and Velani may not seem too restrictive. However, as well as giving a cleaner statement, Theorem \ref{mtp for linear forms} also gives rise to a number of interesting and very general applications in Diophantine Approximation and it is plausible that some of these may have been out of reach when trying to use \cite[Theorem~4.1]{BV slicing}. Applications of Theorem \ref{mtp for linear forms} are discussed at some length in \cite{Allen-Beresnevich} and \cite{Allen-Troscheit}. Further recent applications of Theorem \ref{mtp for linear forms} have also come to light in \cite{Allen-Ramirez,HKS,Zhang}.

\section{Mass Transference Principle with ``Local Scaling Property''}

In \cite{Allen-Baker}, a general mass transference principle was established which emcompasses the aforementioned Theorems \ref{MTP}, \ref{BV general MTP}, and \ref{mtp for linear forms} and further extends these results in several directions. In particular, the general mass transference principle established in \cite{Allen-Baker} holds in locally compact metric spaces, enables the transference of Hausdorff $g$-measure statements to Hausdorff $f$-measure statements, and is applicable to a wider class of sets than are covered by the previously mentioned theorems. For example, the result established in \cite{Allen-Baker} (Theorem \ref{Allen-Baker MTP} below) holds for limsup sets defined by smooth manifolds in $\R^k$ as well as self-similar sets arising from iterated function systems satisfying the open set condition.

While Theorems \ref{MTP} and \ref{BV general MTP} apply to limsup sets of balls and Theorem \ref{mtp for linear forms} deals with essentially limsup sets of neighbourhoods of planes, the result in \cite{Allen-Baker} is applicable to limsup sets defined by neighbourhoods of sets satisfying a certain \emph{local scaling property}. Essentially, this local scaling property that is introduced ensures that the sets under consideration behave locally like affine subspaces. This property appears to be the key feature enabling one all-encompassing statement which includes and extends several of the aforementioned results.

Let $(X,d)$ be a locally compact metric space and let $g$ be a doubling dimension function which furthermore satisfies \eqref{g doubling condition}. Given a set $F \subset X$ and a real number $\delta \geq 0$, we define the $\delta$-neighbourhood of $F$ to be 
\[\Delta(F, \delta) = \{x \in X: d(x,F)<\delta\},\]
where $d(x,F) = \min\{d(x,y): y \in F\}$.

\begin{definition}[Local Scaling Property] \label{LSP definition}
Let $\cF=(F_j)_{j \in \N}$ be a sequence of sets in $X$ and suppose $0 \leq \kappa < 1$ is a real number. We say that $\cF$ satisfies the \emph{local scaling property} (henceforth abbreviated to \emph{LSP}) with respect to $\kappa$ if there exist constants $c_3,c_4>0$ and $r_0>0$ such that, for any real numbers $r$ and $\delta$ with $0 < \delta < r < r_0$, and for any $j \in \N$ and $x \in F_j$, we have that 
\begin{align} \label{LSP}
c_3 g(\delta)^{1-\kappa} g(r)^{\kappa} \leq \cH^g(B(x,r) \cap \Delta(F_j,\delta)) \leq c_4 g(\delta)^{1-\kappa} g(r)^{\kappa}.
\end{align}
If \eqref{LSP} is satisfied for one particular set $F \subset X$, we will also say that $F$ satisfies the LSP with respect to $\kappa$. It should be clear from context when we are referring to a single set satisfying the LSP and when we are referring to a collection of sets satisfying the LSP.
\end{definition}

The following examples of sets satisfying the LSP are given in \cite{Allen-Baker}. We merely state the examples here and, in the case of Examples \ref{LSP example: manifolds} and \ref{LSP example: IFSs}, refer the reader to \cite{Allen-Baker} for more detailed arguments for why such sets satisfy the LSP.

\begin{example}
Suppose $\cF$ is a sequence of points in some locally compact metric space $(X,d)$. Then $\cF$ satisfies the LSP with respect to $\kappa = 0$. 
\end{example}

\begin{example}
Suppose that $F \subset \R^k$ is any set satisfying the LSP with respect to some $0 \leq \kappa < 1$ and suppose that $\cF = (F_j)_{j \in \N}$ is a sequence of sets where, for each $j \in \N$, $F_j$ is the image of $F$ under some isometry. Then $\cF$ satisfies the LSP with respect to $\kappa$.   
\end{example}

\begin{example} \label{LSP example: manifolds}
Let $\cM \subset \R^k$ be a smooth compact $\ell$-dimensional manifold. Then $\cM$ satisfies the LSP with respect to $\kappa = \frac{\ell}{k}$.     
\end{example}

\begin{example} \label{LSP example: IFSs}
Suppose $N \in \N$ and that $\Phi = \{\varphi_j:\R^k \to \R^k\}_{j=1}^{n}$ is a self-similar iterated function system. Denote by $K$ the self-similar set associated with $\Phi$. Furthermore, suppose that $\Phi$ satisfies the open set condition. Then $K$ satisfies the LSP with respect to $\kappa = \frac{\dimh{K}}{k}$.
\end{example}

Now we turn to stating the main result of \cite{Allen-Baker}. Suppose that $\cF = (F_j)_{j \in \N}$ is a sequence of sets in $X$ and $\Upsilon = (\Upsilon_j)_{j \in \N}$ is a sequence of non-negative real numbers with $\Upsilon_j \to 0$ as $j \to \infty$. We define the limsup set
\begin{align*}
\Lambda(\Upsilon) &= \{x \in X: x \in \Delta(F_j,\Upsilon_j) \quad \text{for infinitely many } j \in \N\}.
\end{align*}
The following general mass transference principle was established in \cite[Theorem 1]{Allen-Baker}.

\begin{theorem} \label{Allen-Baker MTP}
Let $(X,d)$ be a locally compact metric space and let $g$ be a doubling dimension function which also satisfies \eqref{g doubling condition}. Let $\cF=(F_j)_{j \in \N}$ be a sequence of sets in $X$ satifying the LSP with respect to some $0 \leq \kappa < 1$, and let $\Upsilon = (\Upsilon_j)_{j \in \N}$ be a sequence of non-negative real numbers such that $\Upsilon_j \to 0$ as $j \to \infty$. Suppose that $f$ is a dimension function such that $f/g$ is monotonic and $f/g^{\kappa}$ is also a dimension function. 

Suppose that for any ball $B$ in $X$ we have
\begin{align*}
\cH^g\left(B \cap \Lambda\left(g^{-1}\left(\left(\frac{f(\Upsilon)}{g(\Upsilon)^{\kappa}}\right)^{\frac{1}{1-\kappa}}\right)\right)\right) = \cH^g(B).
\end{align*}
Then, for any ball $B$ in $X$, we have 
\[\cH^f(B \cap \Lambda(\Upsilon)) = \cH^f(B).\]
\end{theorem}

The strategy used for proving Theorem \ref{Allen-Baker MTP} mirrors the ideas used in \cite{Allen-Beresnevich} for proving Theorem \ref{mtp for linear forms}. In turn, the ideas there are inspired by those used by Beresnevich and Velani in~\cite{BV MTP} for proving the ``original'' Mass Transference Principle (Theorem \ref{MTP}). If we take $\cF$ to be a sequence of points in a locally compact metric space and $\kappa = 0$, then Theorem \ref{Allen-Baker MTP} matches Theorem \ref{BV general MTP}. If our underlying metric space is $\R^k$, $g(r)=r^k$, $\cF$ is a sequence of points in $\R^k$, and $\kappa=0$, then Theorem \ref{Allen-Baker MTP} reduces to Theorem \ref{MTP}. We can also derive Theorem \ref{mtp for linear forms} from Theorem \ref{Allen-Baker MTP} by taking our metric space to be a suitable ball $\Omega$ in $\R^k$, $g(r)=r^j$, $\cF$ to be a sequence of $\ell$-dimensional planes in $\R^k$, and $\kappa = \frac{\ell}{k}$. Thus, Theorem \ref{Allen-Baker MTP} incorporates several of the previously discovered mass transference principles. However, as can be seen from Examples \ref{LSP example: manifolds} and \ref{LSP example: IFSs} above, Theorem \ref{Allen-Baker MTP} is also more widely applicable to much more general sets.

\section{Mass Transference Principle from Balls to Rectangles} \label{balls to rectangles section}

The variants of the mass transference principle we exhibit in this section and the next allow one to derive Hausdorff measure results for limsup sets arising from rectangles. These results are particularly relevant to the topic of weighted approximation in various contexts and have already found several applications in this direction \cite{Allen-Wang, BLW2021, BLW2023, HW, WL, WLWY, Wu}. 

In this section we will present the mass transference principle from ``balls to rectangles'' proved by Wang, Wu, and Xu in \cite{WWX MTP}. In the next section we will discuss mass transference principles from ``rectangles to rectangles''.

Let $k \in \N$ and suppose that $(x_n)_{n \in \N}$ is a sequence of points in the $k$-dimensional unit cube~$[0,1]^k$. Let $(r_n)_{n \in \N}$ be a sequence of non-negative real numbers such that $r_n \to 0$ as $n \to \infty$. We define
\begin{align*}
    W_{0} &= [0,1]^k \cap \limsup_{n \to \infty}{B(x_n, r_n)} \\ 
          &= \left\{x \in [0,1]^k: x \in B(x_n,r_n) \text{ for i.m. } n \in \N\right\}.    
\end{align*}
For $\ba = (a_1,a_2,\dots,a_k) \in \R^{k}_{\geq 0}$, $x = (x_1,\dots,x_k) \in [0,1]^k$, and $r>0$, define
\[B^{\ba}(x,r) = \prod_{i=1}^{k}{B(x_i,r^{a_i})}. \]
In other words, $B^{\ba}(x,r)$ is the rectangle with centre $x$ and side-lengths $(r^{a_1},\dots,r^{a_k})$. Given $(x_n)_{n \in \N}$ and $(r_n)_{n \in \N}$ as above and $\ba = (a_1,\dots,a_k) \in \R^k_{\geq 0}$ with $1 \leq a_1 \leq a_2 \leq \dots \leq a_k$, we define 
\begin{align*}
    W_{\ba} &= [0,1]^k \cap \limsup_{n \to \infty}{B^{\ba}(x_n, r_n)} \\ 
            &= \left\{x \in [0,1]^k: x \in B^{\ba}(x_n,r_n) \text{ for i.m. } n \in \N\right\}.
\end{align*}

In \cite[Theorem 1.2]{WWX MTP}, Wang, Wu, and Xu proved the following result regarding Hausdorff dimension.

\begin{theorem}
Let $(x_n)_{n \in \N}$ be a sequence of points in the $k$-dimensional unit cube $[0,1]^k$, let $(r_n)_{n \in \N}$ be a sequence of non-negative real numbers such that $r_n \to 0$ as $n \to \infty$, and let  $\ba = (a_1,\dots,a_k) \in \R^k_{\geq 0}$ be such that $1 \leq a_1 \leq a_2 \leq \dots \leq a_k$. Suppose that $\cL(W_0) = 1$. Then
\[\dimh{W_{\ba}} \geq \min_{1 \leq j \leq k}\left\{\frac{k + ja_j - \sum_{i=1}^{j}{a_i}}{a_j}\right\}.\]
\end{theorem}

Subject to the additional constraint that $a_k > 1$, Wang, Wu and Xu proved the following result regarding Hausdorff measure in \cite[Theorem 1.3]{WWX MTP}. 

\begin{theorem}
Let $(x_n)_{n \in \N}$ be a sequence of points in the $k$-dimensional unit cube $[0,1]^k$, let $(r_n)_{n \in \N}$ be a sequence of non-negative real numbers such that $r_n \to 0$ as $n \to \infty$, and let  $\ba = (a_1,\dots,a_k) \in \R^k_{\geq 0}$ be such that $1 \leq a_1 \leq a_2 \leq \dots \leq a_k$ and $a_k >1$. Suppose that $\cL(W_0) = 1$ and write 
\[s = \min_{1 \leq j \leq k}\left\{\frac{k + ja_j - \sum_{i=1}^{j}{a_i}}{a_j}\right\}.\]
Then $\cH^s{(W_{\ba})}=\infty$.
\end{theorem}

\section{Mass Transference Principles from Rectangles to Rectangles} \label{rectangles to rectangles section}

The results from the previous section, giving us mass transference principles from ``balls to rectangles'' have been further generalised recently by Wang and Wu in \cite{WW2021} allowing for mass transference principles from ``rectangles to rectangles''.

\subsection{The setup} \label{rectangles setup}

Fix $k \in \N$ and, for each $1 \leq i \leq k$, let $(X_i, |\cdot|_i, m_i)$ be a bounded locally compact metric measure space. Further suppose, for each $1 \leq i \leq k$, that $m_i$ is a $\delta_i$-Ahlfors-David regular probability measure for some $\delta_i >0$. In \cite{WW2021}, Wang and Wu establish mass transference principles from rectangles to rectangles in the product space $(X, |\cdot|, m)$, where
\begin{align*}
    X = \prod_{i=1}^{k}{X_i}, \qquad m = \prod_{i=1}^{k}{m_i}, \qquad \text{and} \quad |\cdot| = \max_{1 \leq i \leq k}{|\cdot|_i}.
\end{align*}

We note that a ball $B(x,r)$ in $X$, where $x= (x_1,\dots,x_k) \in X$ and $r > 0$, is a product of balls from the spaces $\{X_i\}_{i=1}^{k}$, i.e.
\[B(x,r) = \prod_{i=1}^{k}{B(x_i,r)}.\]

Next, let $J$ be an infinite countable indexing set and let $\beta: J \to \R_{>0}$ be such that for any real number $M > 1$ we have
\[\#\{\alpha \in J: \beta(\alpha)<M\} < \infty.\] 
Let $(\ell_n)_{n \in \N}$ and $(u_n)_{n \in \N}$ be sequences of integers such that $u_n \geq \ell_n$ for all $n \in \N$ and $\ell_n \to \infty$ as $n \to \infty$. One implication of these assumptions will be that without loss of generality we may assume that $\ell_n > 0$ and $u_n > 0$ for all $n \in \N$. For $n \in \N$, define
\[J_n = \{\alpha \in J: \ell_n \leq \beta(\alpha) \leq u_n\}.\]

For each $1 \leq i \leq k$, let $(\cR_{\alpha,i})_{\alpha \in J}$ be a collection of subsets in $X_i$. We define the corresponding sequence of \emph{resonant sets} in $X$ by setting 
\[\cR_{\alpha} = \prod_{i=1}^{k}{\cR_{\alpha,i}}.\]
Let $\rho: \R_{>0} \to \R_{>0}$ be a non-increasing function such that $\rho(u) \to 0$ as $u \to \infty$. For $\ba = (a_1,\dots,a_k) \in \R_{\geq 0}^{k}$ and $\rho \in \R_{>0}$, define
\[\Delta(\cR_{\alpha}, \rho^{\ba}) = \prod_{i=1}^{k}{\Delta(\cR_{\alpha,i},\rho^{a_i})},\]
where $\Delta(\cR_{\alpha,i},\rho^{a_i})$ denotes the $\rho^{a_i}$-neighbourhood of $\cR_{\alpha,i}$ in $X_i$. Note that $\Delta(\cR_{\alpha}, \rho^{\ba})$ defines a ``rectangle'' in $X$. To establish the mass transference principles given in \cite{WW2021}, Wang and Wu further require that, for each $1 \leq i \leq k$, the collections $(\cR_{\alpha,i})_{\alpha \in J}$ satisfy the local scaling property (see Definition \ref{LSP definition}) with respect to some $0 \leq \kappa < 1$ and where $g(r)=r^{\delta_i}$.

Suppose that $\ba = (a_1,\dots,a_k) \in  \R_{\geq 0}^{k}$ is fixed. Given $\bt = (t_1,\dots, t_k)\in  \R_{\geq 0}^{k}$, define the limsup set 
\[W(\bt) = \left\{x \in X: x \in \Delta(\cR_{\alpha}, \rho(\beta(\alpha))^{\ba + \bt}) \quad \text{for i.m. } \alpha \in J\right\}.\]
This is essentially the limsup set obtained by ``shrinking'' the rectangles $\Delta(\cR_{\alpha},\rho(\beta(\alpha))^{\ba})$. In \cite{WW2021}, Wang and Wu also consider variations on the set $W(\bt)$ with more general approximating functions. However, for simplicity, here we will restrict our attention to the setting outlined above.

In \cite{WW2021}, Wang and Wu essentially provide two ``types'' of mass transference principle from rectangles to rectangles. The first type requires the original sequence of ``unshrunk'' rectangles to satisfy certain \emph{ubiquity} conditions, as defined in the following section. The second type is more reflective of the typical statements of mass transference principles seen in the previous sections. More precisely, the second type of mass transference principle from rectangles to rectangles exhibited by Wang and Wu in \cite{WW2021} allows us to infer information about the Hausdorff dimension of sets of the form $W(\bt)$ (as defined above) subject to knowing that the limsup set of ``orginal'' rectangles is of full measure.

\subsection{Ubiquity for rectangles}

Inspired by the notions of local and global ubiquitous systems defined by Beresnevich, Dickinson, and Velani in \cite{BDV}, Wang and Wu defined the natural corresponding notions of \emph{local ubiquity} and \emph{uniform local ubiquity} for rectangles in \cite{WW2021}. Throughout we assume the same setting and use the same notation as described in the previous subsection.

\begin{definition}[Local ubiquity for rectangles]
We say that $(\{\cR_{\alpha}\}_{\alpha \in J}, \beta)$ is a \emph{locally ubiquitous system of rectangles} with respect to $(\rho, \ba)$ if there exists a constant $c > 0$ such that for any ball $B$ in $X$ we have
\[\limsup_{n \to \infty}{m\left(B \cap \bigcup_{\alpha \in J_n}{\Delta(\cR_{\alpha}, \rho(u_n)^{\ba})}\right)} \geq c m(B). \]
\end{definition}

\begin{definition}[Uniform local ubiquity for rectangles]
We say that $(\{\cR_{\alpha}\}_{\alpha \in J}, \beta)$ is a \emph{uniformly locally ubiquitous system of rectangles} with respect to $(\rho, \ba)$ if there exists a constant $c > 0$ such that for any ball $B$ in $X$, there is some $n_0(B) \in \N$ such that
\[m\left(B \cap \bigcup_{\alpha \in J_n}{\Delta(\cR_{\alpha}, \rho(u_n)^{\ba})}\right) \geq c m(B) \]
for all $n \geq n_0(B)$.
\end{definition}

\subsection{Mass transference principles from rectangles to rectangles}

Throughout this subsection, let us again assume the setting described in Section \ref{rectangles setup}. Recall that for each $1 \leq i \leq k$, we assume that the measure $m_i$ is $\delta_i$-Ahlfors-David regular for some $\delta_i >0$. Also recall that we assume that for each $1 \leq i \leq k$, the collections $(\cR_{\alpha,i})_{\alpha \in J}$ satisfy the local scaling property (Definition \ref{LSP definition}) with respect to some $0 \leq \kappa < 1$ with $g(r)=r^{\delta_i}$. Assuming local ubiquity for rectangles, Wang and Wu first establish the following mass transference principle from rectangles to rectangles in \cite{WW2021}. 

\begin{theorem} \label{WW ubiquity dimension result}
Assume the setting as described above and further suppose that $(\{\cR_{\alpha}\}_{\alpha \in J}, \beta)$ is a locally ubiquitous system of rectangles with respect to $(\rho, \ba)$. Then
\[\dimh{W(\bt)} \geq \min_{A \in \cA}\left\{\sum_{j \in \cK_1}{\delta_j}+\sum_{j \in \cK_2}{\delta_j}+\kappa \sum_{j \in \cK_3}{\delta_j}+(1-\kappa)\frac{\sum_{j \in \cK_3}{a_j \delta_j}-\sum_{j \in \cK_2}{t_j \delta_j}}{A}\right\}\]
where
\[\cA = \{a_i, a_i+t_i: 1 \leq i \leq k\}\]
and, for each $A \in \cA$, the sets $\cK_1, \cK_2, \cK_3$ are defined by
\begin{align*}
\cK_1 &= \{1 \leq j \leq k: a_j \geq A\}, \\
\cK_2 &= \{1 \leq j \leq k: a_j+t_j \leq A\}\setminus \cK_1, \\ 
\cK_3 &=\{1,\dots,k\}\setminus (\cK_1 \cup \cK_2).
\end{align*}
Note that the sets $\cK_1, \cK_2, \cK_3$ give a partition of $\{1,\dots,k\}$.
\end{theorem}

Furthermore, Wang and Wu also obtain the following Hausdorff measure result in \cite{WW2021}. 

\begin{theorem} \label{WW ubiquity measure result}
Under the same assumptions as in Theorem \ref{WW ubiquity dimension result}, for any ball $B$ in $X$, we have
\[\cH^s(B \cap W(\bt)) = \cH^s(B)\]
for
\[s= \min_{A \in \cA}\left\{\sum_{j \in \cK_1}{\delta_j}+\sum_{j \in \cK_2}{\delta_j}+\kappa \sum_{j \in \cK_3}{\delta_j}+(1-\kappa)\frac{\sum_{j \in \cK_3}{a_j \delta_j}-\sum_{j \in \cK_2}{t_j \delta_j}}{A}\right\}.\]
\end{theorem}

Assuming uniform local ubiquity for rectangles, Wang and Wu generalise Theorems \ref{WW ubiquity dimension result} and \ref{WW ubiquity measure result} in \cite{WW2021} to 
analogous results with more general approximating functions. However, we shall not elaborate further on these more general statements here.

We conclude this section by stating the ``full measure'' mass transference principle from rectangles to rectangles also established by Wang and Wu in \cite{WW2021}. 

\begin{theorem} \label{WW full measure mtp}
Assume the setting as described above and further suppose that
\[m\left(\limsup_{\substack{\alpha \in J \\ \beta(\alpha) \to \infty}}{\Delta(\cR_{\alpha}, \rho(\beta(\alpha))^{\ba})}\right) = m(X).\]
Then, 
\[\dimh{W(\bt)} \geq \min_{A \in \cA}\left\{\sum_{j \in \cK_1}{\delta_j}+\sum_{j \in \cK_2}{\delta_j}+\kappa \sum_{j \in \cK_3}{\delta_j}+(1-\kappa)\frac{\sum_{j \in \cK_3}{a_j \delta_j}-\sum_{j \in \cK_2}{t_j \delta_j}}{A}\right\}\]
where
\[\cA = \{a_i, a_i+t_i: 1 \leq i \leq k\}\]
and, for each $A \in \cA$, the sets $\cK_1, \cK_2, \cK_3$ are defined by
\begin{align*}
\cK_1 &= \{1 \leq j \leq k: a_j \geq A\}, \\
\cK_2 &= \{1 \leq j \leq k: a_j+t_j \leq A\}\setminus \cK_1, \\ 
\cK_3 &=\{1,\dots,k\}\setminus (\cK_1 \cup \cK_2).
\end{align*}
\end{theorem}
\noindent Note that the conclusion in Theorem \ref{WW full measure mtp} is the same as that in Theorem \ref{WW ubiquity dimension result}.

\section{Mass Transference Principles from Balls to Open Sets}

In this section we present mass transference principles which further generalise some of the theorems seen in earlier sections. In particular, we present results which allow us to pass from Lebesgue measure statements for limsup sets defined by sequences of balls to Hausdorff measure and dimension statements for related limsup sets generated by arbitrary open sets.

Fix $k \in \N$, let $(x_n)_{n \in \N}$ be a sequence of points in $[0,1]^k$ and let $(r_n)_{n \in \N}$ be a sequence of positive real numbers such that $r_n \to 0$ as $n \to \infty$. Consider the sequence of balls $(B_n:=B(x_n ,r_n))_{n\in\mathbb{N}}$ in $[0,1]^k$ and suppose that 
\[\mathcal{L}\left(\limsup_{n \to \infty}{B_n}\right)=1.\] 
In this situation, it follows from the mass transference principles presented earlier in Sections \ref{MTP section} and \ref{balls to rectangles section} that we can estimate 
\[\dimh \left(\limsup_{n \to \infty}{B(x_n ,r_n^{\delta})}\right)\]
where $\delta > 0$ and
\[\dimh \left(\limsup_{n\to\infty}(B^{\ba}(x_n,r_n))\right)\] 
where $\ba = (a_1, \dots, a_k)\in \R^k_{\geq 0}$ and $1 \leq a_1 \leq ... \leq a_k$. The notation $B^{\ba}(x,r)$ here is as defined in Section \ref{balls to rectangles section}. It is natural to ask whether one can derive a more general mass transference principle when one shrinks each ball $B_n$ into a set $U_n \subset B_n$ subject to some appropriate geometric constraints. It is easy to see that if the sets $U_n$ have empty interior, then not much can be said in general so we will consider the case where the sets $U_n$ are open. In this case Koivusalo and Rams \cite{KR} obtained a lower bound for 
\[\dimh \left(\limsup_{n \to \infty}{U_n}\right)\] 
in terms of the Hausdorff content of the sets $U_n$. More precisely, they established  the following result in \cite{KR}.

 \begin{theorem}
 \label{lowKR}
Let $(B_n:=B(x_n ,r_n))_{n\in\N}$ be a sequence of balls in $[0,1]^k$ with $r_n \to 0$ as $n \to \infty$ and suppose that 
\[\cL\left(\limsup_{n \to \infty}{B_n}\right) = 1.\]
Let $(U_n)_{n\in\N}$ be sequence of open subsets of $[0,1]^k$ such that for every $n\in\N$, we have $U_n \subset B_n$. Suppose that there exists a real number $0 < s < k$ such that for every $n \in\mathbb{N}$ we have 
\begin{equation}
\mathcal{H}^s_{\infty}(U_n)\geq \mathcal{L} (B_n).    
\end{equation}
Then
\[\dimh \left(\limsup_{n\to\infty} U_n\right) \geq s.\]
 \end{theorem}
 
Koivusalo and Rams in fact prove a different statement in \cite{KR} from which they deduce Theorem \ref{lowKR} as an immediate corollary. The statement which they prove is formulated in terms of a generalised version of the singular value function (defined by \eqref{definitionsingularvalue}). The advantage of the generalised singular value function introduced in \cite{KR} is that it applies to any set rather than being limited to rectangles or ellipsoids. To be more precise, let $E \subset \R^k$ be a Borel set and let $s \geq 0$ be a real number. We define the \emph{generalised singular value function} of $E$ as
\[\phi_s (E)=\sup_{\mu\in\mathcal{M}(E)}\inf_{x\in E}\inf_{r>0}\frac{r^s}{\mu(B(x,r))}.\]
With this definition, Koivusalo and Rams were able to prove the following result in \cite{KR}.

\begin{theorem} \label{KR MTP}
Let $(B_n:=B(x_n ,r_n))_{n\in\N}$ be a sequence of balls in $[0,1]^k$ with $r_n \to 0$ as $n \to \infty$ and suppose that 
\[\cL\left(\limsup_{n \to \infty}{B_n}\right) = 1.\] Let $(U_n)_{n\in\N}$ be sequence of open subsets of $[0,1]^k$ such that for every $n\in\N$, we have $U_n \subset B_n$. Suppose that there exists a real number $s \geq 0$ such that for every $n \in\mathbb{N}$ we have 
\begin{equation}
\phi_s(U_n)\geq \cL(B_n).    
\end{equation}
Then
\[\dimh \left(\limsup_{n\to\infty} U_n\right) \geq s.\]    
\end{theorem}

That Theorem \ref{lowKR} is a corollary of Theorem \ref{KR MTP} follows from the observation of Koivusalo and Rams in \cite[Proposition 2.1]{KR} that Hausdorff content and the singular value function agree up to a multiplicative constant.

Notice that Theorems \ref{lowKR} and \ref{KR MTP} only relate to the Hausdorff dimension of $\limsup_{n \to \infty}U_n$. A natural refinement of Theorem \ref{KR MTP}, determining the Hausdorff $f$-measures of $\limsup_{n\to\infty}U_n$, was established by Zhong in \cite{Zhong}. To obtain Hausdorff $f$-measure results for general dimension functions, Zhong first further extended the definition of the generalised singular value function in the following way. Let $f: \R_{>0} \to \R_{>0}$ be a dimension function, let $E \subset \R^k$ and let $s \geq 0$. We define the \emph{generalised $f$-singular value function} of $E$ as
\[\phi_f (E)=\sup_{\mu\in\mathcal{M}(E)}\inf_{x\in E}\inf_{r>0}\frac{f(r)}{\mu(B(x,r))}.\]
Zhong proved the following in \cite{Zhong}.
 
 \begin{theorem}
Let $(B_n:=B(x_n ,r_n))_{n\in\N}$ be a sequence of balls in $[0,1]^k$ with $r_n \to 0$ as $n \to \infty$ and suppose that 
\[\cL\left(\limsup_{n \to \infty}{B_n}\right) = 1.\] 
Let $(U_n)_{n\in\N}$ be sequence of open subsets of $[0,1]^k$ such that for every $n\in\N$, we have 
\[U_n \subset B_n \qquad \text{and} \qquad \lim_{n \to \infty}{\frac{\cL(B_n)}{\cL(U_n)}}= \infty.\]
Let $f: \R_{>0} \to \R_{>0}$ be a dimension function such that $r^{-k}f(r)$ is monotonically decreasing. Suppose that for every $n \in\mathbb{N}$ we have 
\begin{equation}
\label{condkzhong}
\phi_f(U_n)\geq \cL(B_n).    
\end{equation}
Then, for any ball $B$ in $[0,1]^k$, we have
\[\cH^f\left(B \cap \limsup_{n \to \infty}{U_n}\right) = \cH^f(B).\]    
 \end{theorem}

We add to this section that a comparable result was also obtained by Persson in \cite{Pers}. In order to state this result, let us define the \emph{Riesz energy} of a set $E$.

\begin{definition}
    Let $E\subset \mathbb{R}^k$ be a Borel set and let $0\leq s\leq k$. We define the \emph{$s$-dimensional Riesz energy} of $E$ as 
    $$I_s(E)=\int_{E}\int_{E}\frac{1}{\vert x-y\vert^s}dxdy.$$
\end{definition}
We also need to define the \emph{class of large intersection of dimension $s$}. 
\begin{definition}
\label{sdimlargeinter}
Let $0\leq s\leq k$ be a real number. Then the \emph{$s$-dimensional class of large intersection of $[0,1]^k$}, denoted by $\cG_s ([0,1]^k)$, is defined to be the set of all $G_{\delta}$-sets (countable intersections of dense open sets) $E$ such that
\[\cH^t_{\infty}(E\cap B(x,r))=r^t\]
for every $t<s$ and for every ball $B(x,r)\subset [0,1]^k$.
\end{definition}
Persson proved the following result in \cite{Pers}.

\begin{theorem} \label{MTPperss}
Let $(B_n:=B(x_n ,r_n))_{n\in\N}$ be a sequence of balls in $[0,1]^k$ with $r_n \to 0$ as $n \to \infty$ and suppose that
\[\cL\left(\limsup_{n\to\infty}B_n\right)=1.\] 
Let $(U_n)_{n\in\mathbb{N}}$ be a sequence of open sets in $[0,1]^k$ such that $U_n \subset B_n$ for every $n\in\N$. Let 
    $$s=\sup\left\{t>0 : \sup_{n \in \N}{\frac{I_t(U_n)\cL(B_n)}{\cL(U_n)^2}}<\infty\right\}.$$
    Then 
\[\dimh \limsup_{n\to\infty}U_n \geq s,\]
and
\[\limsup_{n\to\infty}U_n \in \mathcal{G}_s ([0,1]^k).\] 
\end{theorem}
It is worth noting that the bound provided by Theorem \ref{MTPperss} does not improve upon Theorems \ref{lowKR} and \ref{KR MTP} in terms of dimension, see \cite{Pers} for further discussion regarding this point. 
However, Theorem \ref{MTPperss} implies that the set $\limsup_{n\to\infty}U_n$ belongs to $\mathcal{G}_s([0,1]^k)$, which is not established by Theorem \ref{lowKR}. 

In the setting of general metric spaces, the following mass transference principle from balls to arbitrary open sets has recently been established by Eriksson-Bique in \cite{EB1}.
 
\begin{theorem}
\label{BiqueMTP}
Let $(X,d)$ be a complete metric space, let $s\geq 0$ be a real number and let $(B(x_n ,r_n))_{n\in\N}$ be a sequence of balls in $X$ with $r_n \to 0$ as $n \to \infty$. Suppose that, for every ball $B=B(x,r)$ in $X$, we have
\begin{equation}
\label{condBique}
\mathcal{H}^s_{\infty}\left(B \cap \limsup_{n\to\infty}B(x_n,r_n)\right)\geq C r^s.
\end{equation}
Suppose that $(U_n)_{n\in\N}$ is an arbitrary sequence of open sets in $X$ such that $U_n \subset B(x_n, r_n)$ for every $n\in\N$. Furthermore, suppose that there exists some $0\leq t\leq s$ such that for every $n \in \N$ we have
\begin{equation}
\label{condUn}
\cH^t_{\infty}(U_n )\geq r_n^s.
\end{equation}
Then, there exists a constant $\kappa>0$ such that for every ball $B$ in $X$, we have
\[\mathcal{H}^t_{\infty}\Big(B\cap \limsup_{n\to\infty}U_n \Big)\geq \kappa \mathcal{H}^t_{\infty}(B).\]
In particular, $\dimh \limsup_{n\to\infty}U_n \geq t.$
\end{theorem}
  
 Note that, by \eqref{condBique}, the hypothesis on the sequence $(B_n)_{n\in\mathbb{N}}$ is comparable to belonging to the $s$-dimensional class of large intersection. Moreover, by \eqref{condUn}, in the setting of Theorem~\ref{BiqueMTP}, the set $\limsup_{n\to\infty}U_n$ also satisfies a property that is comparable to belonging to the $t$-dimensional class of large intersection.  For the sake of simplicity, Theorem \ref{BiqueMTP} is stated only for traditional Hausdorff $s$-measure but we remark for completeness that Eriksson-Bique obtains results for more general dimension functions in \cite{EB1}.
 
\section{Mass Transference Principles for Inhomogeneous Measures}

In this section we discuss mass transference principles where the measures involved are not necessarily Ahlfors-David regular. We will refer to measures which are not Ahlfors-David regular as \emph{inhomogeneous measures}.

\subsection{Multifractal mass transference principles}
In many dynamical settings and in multifractal analysis, limsup sets generated by balls well distributed according to a multifractal measure arise naturally and various authors have developed tools similar to the mass transference principle of Beresnevich and Velani \cite{BV MTP} in order to deal with such cases. The first result of this kind was established in 2007 by Barral and Seuret in \cite{BS1}. Barral and Seuret obtained dimension estimates similar to those obtained by Jaffard \cite{Jaffard1,Ja} and Beresnevich and Velani \cite{BV MTP}, but where the measure involved is a multinomial measure, a Gibbs measure on the dyadic grid, a multiplicative cascade or, more generally, a measure which carries enough self-similarity. The particular case where the measure is a Gibbs measure on the dyadic grid was also established independently later on in 2013  by Fan, Shmeling and Troubetzkoy in \cite{FST} when studying shrinking targets associated with the doubling map on $\mathbb{T}^1$.   

In order to state the results of \cite{BS1} and \cite{FST}, let us first recall several relevant definitions. Let $f:\left\{0,1\right\}^{\mathbb{N}}\to\mathbb{R}$ be a H\"{o}lder continuous function, let $k\in\N$ and let $\bi=(i_1,...,i_k)\in \left\{0,1\right\}^k.$ Denote by $\sigma$ the \emph{left-shift} on $\left\{0,1\right\}^{\mathbb{N}}$; that is, 
\[\sigma(i_1,i_2,i_3,i_4,\dots) = (i_2,i_3,i_4,i_5,\dots)\]
for $(i_1,i_2,i_3,i_4,\dots) \in \{0,1\}^{\N}$.
For finite words $\bi = (i_1,\dots,i_{|\bi|}) \in \bigcup_{N=0}^{\infty}{\{0,1\}^{N}}$ we define the \emph{cylinder set} determined by $\bi$ as 
\[[\bi]=\left\{\bj \in \{0,1\}^{\N}: i_{\ell}=j_{\ell} \quad \text{for } 1 \leq \ell \leq |\bi|\right\}.\] 
Next, define
\begin{align*}
Y_{k,f}(\bi)=\sup_{\bj\in [\bi]}e^{\sum_{0\leq i\leq k-1}f(\sigma^i (\bj))}
\end{align*}
and
\begin{align*}
Z_{k,f}=\sum_{\bi\in \left\{0,1\right\}^k}Y_{k,f}(\bi).
\end{align*}
We define the pressure of $f$ as 
\begin{equation}
\label{press}
P_f =\lim_{k\to\infty}\frac{\log Z_{k,f}}{k }.
\end{equation}

The \emph{Gibbs measure} associated with the potential $f$ is defined to be the unique probability  measure on $\Lambda^{\mathbb{N}},$ $\nu_f $, such that for every $k\in\mathbb{N}$, for every $\bi\in\left\{0,1\right\}^k$ and every $\bj \in [\bi]$, we have
 \begin{equation}
 \label{equagibbs}
 C_1 e^{\sum_{0\leq i\leq k-1}f(\sigma^i (\bj)-kP_{f}}\leq \nu_{f}([\bi])\leq C_2 e^{\sum_{0\leq i \leq k-1}f(\sigma^i (\bj))-kP_{f}}.
 \end{equation}
 for some constants $C_1 , C_2 >0$. It is proved in \cite{Bowen} that such a measure exists.

For $\bi = (i_1,i_2,i_3,i_4,\dots) \in \left\{0,1\right\}^{\mathbb{N}}$, the \emph{projection associated with the dyadic grid} $\pi: \left\{0,1\right\}^{\mathbb{N}}\to [0,1]$ is defined as 
$$\pi(\bi)= \sum_{n \in \N} 2^{-n}i_n.$$
 
The following result was established by Barral and Seuret in \cite{BS1} with the aforementioned special case being proved independently later by Fan, Schmeling and Troubetzkoy in \cite{FST}.

\begin{theorem}
\label{gibbsmtp}
 Let $f$ be a potential on  $\left\{0,1\right\}^{\mathbb{N}}$, let $\nu_{f}$ the Gibbs measure associated with $f$, and let $\mu_{f}=\nu_{f}\circ \pi^{-1}.$  Let $(B_n:=B(x_n,r_n))_{n\in\mathbb{N}}$  be any sequence of balls (i.e. intervals) in $[0,1]$ such that $r_n \to 0$ as $n \to\infty$ and 
 \[\mu_{f}(\limsup_{n\to\infty}B_n)=1.\]
 Then, for any $\delta\geq 1$, we have
$$ \dimh\left(\limsup_{n\to\infty}B(x_n,r_n ^{\delta})\right)\geq \frac{\dim(\mu_{f})}{\delta}.$$
\end{theorem}

Theorem \ref{gibbsmtp} gives rise to the following ``multifractal mass transference principle'' as observed in \cite{BS1,FST}.

\begin{corollary}
     Let \[\mathcal{D}=\left\{\left[\frac{i}{2^k},\frac{i+1}{2^k}\right]:n\in\mathbb{N},0\leq i\leq 2^{k}-1\right\}\] be the set of closed dyadic intervals in $[0,1]$ and let $\mu\in\mathcal{M}(\mathbb{R})$ a Gibbs measure associated with the dyadic grid. For  \[D = \left[\frac{i}{2^{k_D}}, \frac{i+1}{2^{k_D}}\right]\in\mathcal{D},\] define \[\varepsilon_D =\sqrt{\frac{\log \log k_D}{\log k_D}}.\] 
     For every $\alpha \in\mathbb{R}$ and $\delta\geq 1$, we have
     $$\dimh \limsup_{\substack{D\in\mathcal{D}: \\ 2^{-k_D(\alpha-\varepsilon_D )}\leq \mu(D)\leq 2^{-k_D(\alpha-\varepsilon_D )}}} D^{\delta}=\frac{\tau_{\mu}^*(\alpha)}{\delta}.$$
\end{corollary}

As the multifractal analysis of Gibbs measures is naturally correlated with the multifractal analysis of the Birkhoff average of a given potential, Theorem \ref{gibbsmtp} turns out to be  particularly well suited to the study of shrinking targets which are dynamically defined. For example, Theorem \ref{gibbsmtp} has been used to study shrinking targets associated with the doubling map by Fan, Schmeling and Troubetzkoy in \cite{FST} and shrinking targets associated with expanding Markov maps by Liao and Seuret in~\cite{LS}. We conclude this section by noting that other multifractal mass transference principle type statements have been obtained in the random (rather than deterministic) setting by Ekstr\"{o}m and Persson \cite{EP} and Persson and Rams \cite{PR}.

\subsection{A  mass transference principle subject to a full inhomogeneous measure statement}

When dealing with inhomogeneous measures (i.e, measures which are not Ahlfors-David regular), one needs to introduce a geometric quantity which plays the same role as the Hausdorff content does in Theorem \ref{lowKR}. The following quantity was introduced in \cite{ED1} for this purpose. Let $A\subset \mathbb{R}^k$ be a Borel set and let $\mu\in\mathcal{M}(\mathbb{R}^k).$ For $s\geq 0,$ we define the \emph{$\mu$-essential Hausdorff content} of $A$ by setting 
    $$ \mathcal{H}^{\mu,s}_{\infty}(A)=\inf\left\{\mathcal{H}^s_{\infty}(E): E\subset A \text{ and }\mu(E)=\mu(A)\right\}.$$
Using this quantity, the following Theorem was established in \cite{ED1}.
\begin{theorem}
\label{IMTP}
    Let $\mu\in\mathcal{M}(\mathbb{R}^k)$, and suppose $(B_n:=B(x_n ,r_n))_{n\in\mathbb{N}}$ is a sequence of closed balls in $\mathbb{R}^k$ with $r_n \to 0$ as $n \to \infty$ and such that 
    \[\mu\left(\limsup_{n\to\infty}\frac{1}{2}B_n\right)=1.\] 
    Let $(U_n)_{n\in\mathbb{N}}$ be a sequence of open sets such that $U_n \subset B_n$ for every $n\in\mathbb{N}$. Let $0\leq s\leq \underline{\dim} (\mu)$ be a real number such that for every $n\in\mathbb{N}$ we have $$\mathcal{H}^{\mu,s}_{\infty}(U_n)\geq \mu(B_n).$$
    Then there exists a dimension function $f:\R_{>0} \to \R_{>0}$ such that 
    \[\lim_{r\to 0}\frac{\log f(r)}{\log r}=s \quad \text{and} \quad \mathcal{H}^{f}\left(\limsup_{n\to\infty}U_n\right)>0.\]
    In particular, 
    \[\dimh\left( \limsup_{n\to\infty}U_n\right) \geq s.\]
\end{theorem} 
In \cite{ED1}, precise estimates of the $\mu$-essential content are provided. In particular, the following corollary is established in \cite{ED1}.
 \begin{corollary}
 \label{CorollaryED}
 Let $m\geq 2$ be an integer, let $\Phi=\left\{\varphi_1,...,\varphi_m\right\}$ be a self-similar IFS and let $X$ be the attractor of $\Phi$. Let $\mu\in\mathcal{M}(\mathbb{R}^k)$ be a self-similar measure associated with $\Phi$. Recall that $\supp{\mu} = X$. Let $(B_n :=B(x_n ,r_n))_{n\in\mathbb{N}}$ be a sequence of balls in $\R^k$ centered on $X$ such that $r_n\to 0$ and $n \to \infty$ and \[\mu\left(\limsup_{n\to\infty}B_n\right)=1.\] Then for every $\delta\geq 1,$ 
 $$ \dimh \left(\limsup_{n\to\infty}B(x_n, r_n^{\delta})\right)\geq \frac{\dim (\mu)}{\delta}.$$
 \end{corollary}
 Note that this corollary holds without any separation assumption on the iterated function system $\Phi$. Corollary \ref{CorollaryED} has been used in \cite{ED3} to study the dimension of certain shrinking targets associated with some overlapping self-similar iterated function systems.

\paragraph{Acknowledgements.} We thank the organisers of the Fractals and Related Fields IV conference in Porquerolles in September 2022 for an excellent conference on a beautiful island.

\bibliographystyle{plain}

\end{document}